\documentclass{amsart}
\usepackage{amsmath, amscd, amssymb, amsthm}
\usepackage{bbm}
\usepackage{latexsym}
\usepackage{amsfonts}
\usepackage{graphicx}
\usepackage[all,cmtip]{xy}
\usepackage[colorlinks,linkcolor=blue,breaklinks=blue,urlcolor=blue,citecolor=blue,anchorcolor=blue,pagebackref]{hyperref}%
\usepackage{geometry}
\newtheorem{theorem}{Theorem}
\newtheorem{lemma}{Lemma}
\newtheorem{corollary}[theorem]{Corollary}

\newtheorem{conjecture}{Conjecture}

\renewcommand*\backref[1]{}
\renewcommand*\backrefalt[4]{ \ifcase #1 \or (cited on page #2) \else (cited on pages #2) \fi}

\newcommand{\be}{\begin{equation}}
\newcommand{\ee}{\end{equation}}
\newcommand{\bea}{\begin{eqnarray}}
\newcommand{\eea}{\end{eqnarray}}

\newcommand{\vs}{\vspace{0.5cm}}

\def\XXint#1#2#3{{\setbox0=\hbox{$#1{#2#3}{\int}$ }
\vcenter{\hbox{$#2#3$ }}\kern-.6\wd0}}

\begin{document}

\title[BTP metrics with constant holomorphic sectional curvature]{Bismut torsion parallel metrics with constant holomorphic sectional curvature}

\author{Shuwen Chen}
\address{Shuwen Chen. School of Mathematical Sciences, Chongqing Normal University, Chongqing 401331, China}
\email{{3153017458@qq.com}}\thanks{Zheng is the corresponding author. He is partially supported by National Natural Science Foundations of China
with the grant No.12071050 and  12141101, and is supported by the 111 Project D21024.}

\author{Fangyang Zheng}
\address{Fangyang Zheng. School of Mathematical Sciences, Chongqing Normal University, Chongqing 401331, China}
\email{20190045@cqnu.edu.cn; franciszheng@yahoo.com} \thanks{}

\subjclass[2020]{53C55 (primary), 53C05 (secondary)}
\keywords{holomorphic sectional curvature, Chern connection, Bismut connection, Bismut torsion parallel, Vaisman manifolds.}

\begin{abstract}
An old conjecture in non-K\"ahler geometry states that, if a compact Hermitian manifold has constant holomorphic sectional curvature, then the metric must be K\"ahler (when the constant is non-zero) or Chern flat (when the constant is zero). It is known to be true in complex dimension $2$ by the work of Balas and Gauduchon in 1985 (when the constant is negative or zero) and Apostolov, Davidov and Muskarov in 1996 (when the constant is positive). In dimension $3$ or higher, the conjecture is only known in some special cases, such as the locally conformally K\"ahler case (when the constant is negative or zero) by the work of Chen, Chen and Nie, or for complex nilmanifolds with nilpotent $J$ by the work of Li and the second named author.

In this note, we confirm the above conjecture for all non-balanced Bismut torsion parallel (BTP) manifolds. Here the BTP condition means that the  Bismut connection has parallel torsion. In particular, the conjecture is valid for all Vaisman manifolds.
\end{abstract}

\maketitle

\tableofcontents

\markleft{Shuwen Chen and Fangyang Zheng}
\markright{BTP metrics with constant holomorphic sectional curvature}

\section{Introduction and statement of results}\label{intro}

Given a compact Hermitian manifold $(M^n,g)$, denote by $\nabla$ its Chern connection and by $T$, $R$ the torsion and curvature of $\nabla$. The holomorphic sectional curvature of $\nabla$ is defined by
$$ H(X) = \frac{R_{X\bar{X}X\bar{X}}} {|X|^4} $$
where $X$ is  any non-zero complex tangent vector of $(1,0)$ type. As a natural inquiry about Hermitian space forms, the following is a  long standing open question in  non-K\"ahler geometry:

\begin{conjecture}[{\bf Constant Holomorphic Sectional Conjecture}] \label{conj1}
Let $(M^n,g)$ be a compact Hermitian manifold of complex dimension $n\geq 2$. Assume that $H=c$ is a constant. Then $g$ must be K\"ahler when  $c\neq 0$ and  $g$ must be Chern flat (namely, $R=0$) when $c=0$.
\end{conjecture}

In other words, the conjecture says that any compact Hermitian manifold with constant holomorphic sectional curvature must be either K\"ahler (hence a complex space form) or Chern flat. Note that the compactness assumption in the above conjecture is necessary, as there are counterexamples in the non-compact case. 

By the famous result of Boothby \cite{Boothby} in 1958, compact Chern flat manifolds are exactly compact quotients of complex Lie groups (equipped with left-invariant metrics compatible with the complex structure). Such manifolds can (and often) be non-K\"ahler when $n\geq 3$. 

The above conjecture is known to be true when $n=2$, by the work of Balas and Gauduchon \cite{BG} (see also \cite{Balas}) in 1985 for the $c\leq 0$ case, and by 
Apostolov, Davidov, Muskarov \cite{ADM} in 1996 for all $n=2$ cases, as a corollary of their beautiful classification theorem for compact self-dual Hermitian surfaces. 

For $n\geq 3$, the conjecture is still largely open, and there are only a few partial results known so far. In \cite{Tang} K. Tang  confirmed the conjecture under the additional assumption that $g$ is Chern K\"ahler-like (namely the Chern curvature $R$ obeys all K\"ahler symmetries). In \cite{CCN}, H. Chen, L. Chen and X. Nie proved the conjecture under the additional assumption that $g$ is locally conformally K\"ahler and $c\leq 0$. In \cite{ZhouZ}, W. Zhou and the second named author proved that any compact balanced threefold with zero {\em real bisectional curvature}, a notion introduced in \cite{XYangZ} which is slightly stronger than $H$, must be Chern flat. Also, in \cite{LZ} and \cite{RZ}, the second named author and collaborators confirmed the conjecture under the additional assumption that either $(M^n,g)$ is a complex nilmanifold with nilpotent complex structure $J$ (in the sense of \cite{CFGU}), or $g$ is Bismut K\"ahler-like (namely, the curvature $R^b$ of the Bismut connection $\nabla^b$ obeys all K\"ahler symmetries).

Conjecture \ref{conj1} has a twin version, when one replaces the Chern connection by the Riemannian (Levi-Civita) connection $\nabla^r$:

\begin{conjecture} \label{conj2}
Let $(M^n,g)$ be a compact Hermitian manifold of complex dimension $n\geq 2$. Assume that the Riemannian holomorphic sectional curvature $H^r=c$ is a constant. Then $g$ must be K\"ahler when  $c\neq 0$ and  $g$ must be Riemannian flat (namely, $R^r=0$) when $c=0$.
\end{conjecture}

The $n=2$, $c\leq 0$ case for Conjecture \ref{conj2} was proved by Sato and Sekigawa \cite{SS} in 1990, and the $c>0$ case was established by \cite{ADM}. Note that by the classic Bieberbach Theorem, given any compact Riemannian flat manifold, a finite (unbranched) cover of it must be a flat torus. So up to a finite cover, a Riemannian flat compact Hermitian manifold means a (integrable) complex structure $J$ on a flat torus $(T^{2n}_{\mathbb R},g)$ with $J$ being compatible with the flat metric $g$. When $n\leq 2$, $g$ must be K\"ahler and we end up with a flat complex torus. But when $n\geq 3$, there are plenty of such $J$ with $g$ non-K\"aher. In fact in this case $J$ is non-K\"ahlerian, meaning that the complex manifold does not admit any K\"ahler metric. For $n=3$, all such $J$ were classified by the second named author and collaborators in \cite{KYZ}, but for $n\geq 4$, the classification/description of all orthogonal complex structures on  flat tori $T^{2n}_{\mathbb R}$ is still an open question.     

The purpose of this note is confirm the above two conjectures for a special class of Hermitian manifolds, the so-called {\em Bismut torsion-parallel} (BTP for brevity) manifolds. 

Recall that the Bismut connection $\nabla^b$ (also called Strominger connection in some literature) of a Hermitian manifold $(M^n,g)$ is the unique connection satisfying $\nabla^bg=0$, $\nabla^bJ=0$,  and with totally skew-symmetric torsion.  See \cite{Bismut}, \cite{Strominger} for its historic origin. The metric is BTP  if $\nabla^bT^b=0$, or equivalently, if $\nabla^bT=0$. Here we denoted by $T$, $T^b$ the torsion of Chern or Bismut connection, respectively. BTP manifolds form a relatively large and interesting class of special Hermitian manifolds. For instance, all Bismut K\"ahler-like (BKL) manifolds are BTP, as proved by \cite{ZZCrelle}. Here BKL means that the curvature of $\nabla^b$ obeys all K\"ahler symmetries. Clearly, all Bismut flat manifolds (\cite{WYZ}) are BKL. As another example, all Vaisman manifolds are BTP, by the result of Andrada and Villacampa \cite{AndradaV}. Recall that Vaisman means locally conformally K\"ahler manifolds whose Lee form is parallel under the Riemannian (Levi-Cvita) connection. 

In dimension $2$, BTP = BKL = Vaisman. Such surfaces were fully classified by Belgun \cite{Belgun} in 2000. But when $n\geq 3$, BKL and Vaisman are disjoint, and their union is a proper subset of the set of all non-balanced BTP manifolds. Also, when $n\geq 3$, there are BTP manifolds that are balanced (and non-K\"ahler). Balanced BTP manifolds form a highly restrictive and interesting set. It contains examples that are Chern flat or Fano. We refer the readers to the preprints \cite{ZhaoZ24, ZhaoZ25} (which are rewriting of the older file \cite{ZhaoZ22}) and the references therein for more discussions on BTP manifolds.

We will say that {\em Conjecture \ref{conj1} (or \ref{conj2}) holds for a particular class of Hermitian manifolds,} if any manifold in that class cannot have constant Chern (or Riemannian) holomorphic sectional curvature unless it is K\"ahler (hence a complex space form) or it is Chern (or Riemannian) flat. The main result of this note is the following

\begin{theorem} \label{thm1}
Conjectures \ref{conj1} and \ref{conj2} hold for all non-balanced BTP manifolds  and all three-dimensional balanced BTP manifolds.
\end{theorem}

Note that by \cite[Theorem 3]{YangZ},  any compact Hermitian manifold that is Chern (or Riemannian) flat (or more generally, Chern K\"ahler-like or Riemannian K\"ahler-like) must be balanced. So for the class `non-balanced BTP manifolds' in Theorem \ref{thm1}, it simply says that such a manifold cannot have constant Chern (or Riemannian) holomorphic sectional curvature.  Also, since all Vaisman manifolds are non-balanced BTP, an immediate consequence to Theorem \ref{thm1} is the following:

\begin{corollary}
Conjectures \ref{conj1} and \ref{conj2} hold for all compact Vaisman manifolds.
\end{corollary}  

Note that Chen-Chen-Nie \cite{CCN} has confirmed the conjectures for all locally conformally K\"ahler manifolds when the constant $c$ is negative or zero. So the result of the corollary is only new in the case when $c>0$.

After Conjectures \ref{conj1} and \ref{conj2}, it is natural to wonder about what happens when the connection is replaced by $\nabla^b$, which is one of the three canonical metric connections that are widely studied. In other words, one would wonder about how {\em Bismut space forms} would be like. In our previous work \cite{ChenZ} we raised the following:

\begin{conjecture} \label{conj3}
Let $(M^n,g)$ be a compact Hermitian manifold with $n\geq 2$. Assume that Bismut holomorphic sectional curvature $H^b$ is a non-zero constant. Then $g$ must be K\"ahler.
\end{conjecture}

Here $H^b(X)=R^b_{X\bar{X}X\bar{X}}/|X|^4$ is the holomorphic sectional curvature of $R^b$, with $R^b$ the curvature of the Bismut connection $\nabla^b$.

Note that here we omitted the $H^b= 0$ case. As proved in \cite{ChenZ}, for $n=2$ the statement $``H^b=0 \, \Longrightarrow \, R^b=0"$ is valid, but for $n\geq 3$, there are examples of compact Hermitian manifolds with $H^b=0$ but with $R^b\not\equiv 0$. In \cite[Lemma 5]{ChenZ} it was shown that any isosceles Hopf manifold equipped with the standard metric will have $H^b=0$. But for $n\geq 3$, it does not have vanishing $R^b$.  These `counterexamples' on the other hand seem to be rather restrictive, so in \cite{ChenZ} we also asked  the following question:

\vspace{0.25cm}

\noindent {\bf Question 1.} \label{question1}
\emph{ For $n\geq 3$, what kind of compact Hermitian manifold $(M^n,g)$ can have $H^b=0$ but  $R^b\not\equiv 0$?}

\vspace{0.25cm}

Beyond Bismut connection, one may also consider {\em Gauduchon connections}, which is the line joining Chern and Bismut connection: for any real number $t$,   $\nabla^{(t)} = (1-t)\nabla + t \nabla^b$ is called the {\em $t$-Gauduchon connection}.   Denote by $R^{(t)}$ its curvature and by $H^{(t)}$ its holomorphic sectional curvature. A beautiful recent result by Lafuente and Stanfield \cite{LS} states that

\begin{theorem}[Lafuente-Stanfield] \label{thm3}
For any $t\neq 0,1$, a compact $t$-Gauduchon flat manifold must be K\"ahler.
\end{theorem}

In other words, other than Chern and Bismut, any other Gauduchon connection cannot be flat unless the metric is K\"ahler. In fact what they proved is a stronger statement, namely, for any $t\neq 0,1$, if a compact Hermitian manifold is {\em $t$-Gauduchon K\"ahler-like} (see for example \cite{Fu-Zhou} or \cite{ZZGKL}), meaning its $R^{(t)}$ obeys all K\"ahler symmetries, then the metric is K\"ahler. 

Analogous to Conjecture \ref{conj1} or \ref{conj2}, H. Chen and X. Nie in \cite{ChenNie} proposed the following generalization to Conjecture \ref{conj1}:

\begin{conjecture}[{\bf Chen-Nie}] \label{conj4}
Let $(M^n,g)$ be a compact Hermitian manifold of complex dimension $n\geq 2$. Assume that $t\neq 0,1,-1$ and the $t$-Gauduchon holomorphic sectional curvature $H^{(t)}=c$ is a constant. Then $g$ must be K\"ahler.
\end{conjecture}

In \cite{ChenNie}, H. Chen and X. Nie  confirmed Conjecture \ref{conj4} in the $n=2$ case. 

The reason to exclude the $t=-1$ case is that, as in the Bismut case, there are compact Hermitian manifolds with  $H^{(-1)}=0$ but $R^{(-1)}\not\equiv 0$ when $n\geq 3$.  So for this connection, the question should be proposed just like the Bismut case. Therefore  in Conjecture \ref{conj3} and Question \ref{question1} we will include $\nabla^{(-1)}$ there.  

Parallel to Theorem \ref{thm1}, we have the following

\begin{theorem} \label{thm4}
Conjecture \ref{conj3} (for both $\nabla^b$ and $\nabla^{(-1)}$) and Conjecture \ref{conj4} hold for all compact non-balanced BTP manifolds and all compact balanced BTP threefolds.  
\end{theorem}

In summary, what we proved in this note is that,  excluding the Bismut connection $\nabla^b$  and its `mirror' $\nabla^{(-1)}$, for all other Gauduchon connections (including Chern connection) and for the Riemannian connection, the Constant Holomorphic Sectional Conjecture is valid for all non-balanced BTP manifolds (which include BKL and  Vaisman manifolds) and all balanced BTP threefolds. For $\nabla^b$ and $\nabla^{(-1)}$, the same is true when the constant $c\neq 0$, while the $c=0$ case has counterexamples when $n\geq 3$.

At present time, we do not know how to deal with balanced BTP manifolds in dimensions $n\geq 4$, and we do not have a good answer to Question \ref{question1}. We intend to explore them in the future.

\vspace{0.3cm}

\section{Preliminaries}

In this section we will set up the notations and collect some basic facts in Hermitian geometry. 
Let $(M^n,g)$ be a Hermitian manifold. Given any linear connection $D$ on $M^n$, its torsion and curvature are given by
$$ T^D(x,y)=D_xy-D_yx-[x,y], \ \ \ \ \ \ R^D_{xy}z=D_xD_yz-D_yD_xz-D_{[x,y]}z, $$
where $x, y, z$ are vector fields on $M$. They obey the first Bianchi identity
\begin{equation} \label{B1}
{\mathfrak S}\{ R_{xy}z - (D_zT^D)(x,y) - T^D(T^D(x,y),z) \} \ = \ 0,
\end{equation}
where ${\mathfrak S}$ stands for cyclic sum in $x,y,z$. Using the metric $g$, we will write
$$ R^D(x,y,z,w)= g( R^D_{xy}z, w)$$
as a covariant $4$-tensor. Clearly, $R^D$ is skew-symmetric with respect to its first two positions. When  $D$ is a metric connection (namely, when $Dg=0$), then $R^D$ is also  skew-symmetric with respect to its last two positions. So $R^D$ becomes a bi-linear form on $\Lambda^2(TM)$, which in general may not be symmetric as $D$ may have torsion. Nonetheless, one could still define the {\em sectional curvature} of $D$ by
$$ K^D(\pi ) = K^D(x\wedge y ) = \frac{R^D(x,y,y,x)}{|x\wedge y|^2} $$
for any $2$-plane $\pi$  spanned by $x$ and $y$ in the tangent space $T_pM$. Here we wrote  $g=\langle \,, \rangle$ and  as usual $|x\wedge y|^2=|x|^2|y|^2 - \langle x,y\rangle ^2$. Clearly, this value  is independent of the choice of the basis $\{ x,y\}$ of $\pi$. Since our metric $g$ is Hermitian, namely, $\langle Jx, Jy\rangle = \langle x, y\rangle$ for any $x,y$, we can define the {\em holomorphic sectional curvature } $H^D$ of $D$ as the restriction of $K^D$ on all $J$-invariant $2$-planes:
$$ H^D(X) = K^D(x\wedge Jx) = \frac{R^D(x,Jx,Jx,x)}{|x|^4} =  \frac{R^D(X,\overline{X},X,\overline{X})}{|X|^4}  ,$$
where $X=x-iJx$ is the complex tangent vector of type $(1,0)$ corresponding to $x$. Here and below we have extended  $\langle \,, \rangle$ and $R^D$ complex multi-linearly.

Let $e = \{ e_1, \ldots , e_n\}$ be a local unitary frame, namely, each $e_i$ is a complex tangent vector field of type $(1,0)$ and $\langle e_i, \overline{e}_j \rangle = \delta_{ij}$ for $1\leq i,j\leq n$. Denote by $\nabla^r$, $\nabla$, and $\nabla^b$ the Riemannian (Levi-Civita), Chern, and Bismut connection of $(M^n,g)$, by $T$, $T^b$ the torsion of $\nabla$ and $\nabla^b$, and by $R^r$, $R$, $R^b$ the corresponding curvature tensor. Under the frame $e$, write
\begin{equation*}
T(e_i,e_k) = \sum_{j=1}^n T^j_{ik} e_j.
\end{equation*}
Note that our $T^j_{ik}$ here is twice as much as the same notation in references such as \cite{YangZ}, \cite{ZhaoZ22}, so some of the coefficients in formula will differ by a factor of $2$. Under the frame $e$, we have (see for example \cite[Lemma 2]{YangZ})
\begin{eqnarray}
&& \gamma e_i  := \nabla^b e_i- \nabla e_i = \sum_{j,k} \{ T^j_{ik}\varphi_k  - \overline{T^i_{jk}} \overline{\varphi}_k \} e_j,   \label{gamma} \\
&& \nabla^r e_i = \nabla e_i + \frac{1}{2}\gamma e_i +  \frac{1}{2}\sum_{j,k} \{  T^k_{ij}\overline{\varphi}_k  \}  \overline{e}_j,  \label{nablaLC}
\end{eqnarray}
where $\varphi = \{ \varphi_1, \ldots , \varphi_n \}$ is the coframe of local $(1,0)$-forms dual to $e$. 

As is well known, $T(e_i, \overline{e}_k)=0$, and the only possibly non-trivial components of Chern curvature $R$ are $R_{i\bar{j}k\bar{\ell}}$. So the first Bianchi identity (\ref{B1}) for the Chern connection $\nabla$ becomes
\begin{eqnarray} 
&&  {\mathfrak S} \{ T^{\ell}_{ij;k} + \sum_r T^r_{ij}T^{\ell}_{rk} \} \ = \ 0, \label{B1C1} \\
&&  - T^{\ell}_{ik;\bar{j}} \ = \ R_{i\bar{j}k\bar{\ell}} - R_{k\bar{j}i\bar{\ell}}, \label{B1C2}
\end{eqnarray}
for any $1\leq i,j,k,\ell \leq n$, where ${\mathfrak S}$ denotes the cyclic sum in $i,j,k$ and subscripts after the semicolon stand for covariant derivatives with respect to $\nabla$. 

For our later proof, we shall need the following formula \cite[Lemma 7, formula (23)]{YangZ} which relates the components of Riemannian curvature $R^r$ and Chern curvature $R$ under any unitary frame $e$:
\begin{equation} \label{eq:curvature1}
R^r_{i\bar{j}k\bar{\ell}} -  R_{i\bar{j}k\bar{\ell}} = \frac{1}{2} (T^{\ell}_{ik;\bar{j}} + \overline{   T^{k}_{j\ell ;\bar{i}}   }) + \frac{1}{4} \sum_r \{ T^r_{ik} \overline{T^r_{j\ell } }   -  T^{\ell}_{ir} \overline{T^k_{jr} } -  T^j_{kr} \overline{T^i_{\ell r } }  \}  ,
\end{equation}
for any $1\leq i,j,k,\ell \leq n$ and again the subscripts after semicolon stand for covariant derivatives with respect to $\nabla$. (Note that in \cite[Lemma 7]{YangZ}, the Riemannian curvature $R^r$ was denoted as $R$, the Chern curvature $R$ was denoted as $R^h$, and our $T^j_{ik}$ was denoted as $2T^j_{ik}$ there).

Recall that {\em Bismut torsion-parallel manifolds} (BTP for brevity) are Hermitian manifolds whose Bismut connection has parallel torsion: $\nabla^bT^b=0$. This is equivalent to $\nabla^bT=0$ as $T^b$ can be expressed by $T$ (and vice versa). For BTP manifolds, we would like to use Bismut derivatives instead of the Chern derivatives, so we need their relation formula below. By using the formula (\ref{gamma}) for the tensor $\gamma = \nabla^b-\nabla$, a direct computation leads to the following:
\begin{equation} \label{eq:derivative}
T^{\ell}_{ik;\bar{j}} - T^{\ell}_{ik,\bar{j}} \ = \  \sum_r \{ - T^r_{ik} \overline{T^r_{j\ell } }   -  T^{\ell}_{kr} \overline{T^i_{jr} } + T^{\ell}_{ir} \overline{T^k_{j r } }  \},  
\end{equation}
for any $1\leq i,j,k,\ell \leq n$, where the subscripts after comma stand for covariant derivatives with respect to the Bismut connection $\nabla^b$.
Plugging (\ref{eq:derivative}) into (\ref{eq:curvature1}), we obtain the following:

\begin{lemma} \label{lemma1}
Given any Hermitian manifold $(M^n,g)$ and under any local unitary frame $e$, the components of Riemannian curvature and Chern curvature are related by
\begin{equation*} \label{eq:curvature2}
R^r_{i\bar{j}k\bar{\ell}} -  R_{i\bar{j}k\bar{\ell}} = \frac{1}{2} (T^{\ell}_{ik,\bar{j}} + \overline{   T^{k}_{j\ell ,\bar{i}}   }) + \frac{1}{4} \sum_r \{ -3 T^r_{ik} \overline{T^r_{j\ell } }   +3  T^{\ell}_{ir} \overline{T^k_{jr} } -  T^j_{kr} \overline{T^i_{\ell r } } -2 T^{j}_{ir} \overline{T^k_{\ell r} } -2  T^{\ell}_{kr} \overline{T^i_{jr} } \} , 
\end{equation*}
where subscripts after comma stand for covariant derivatives with respect to the Bismut connection $\nabla^b$.
\end{lemma}

Similarly, either by a straight-forward computation using the structure equations and formula (\ref{gamma}) as in  (the $t=0$ case of) Lemma \ref{lemma4} later, or by directly applying \cite[Lemma 3.1]{ZhaoZ22} (and noticing again the factor of $2$ resulted from the notation discrepancy in $T^j_{ik}$), we get the following

\begin{lemma} \label{lemma2}
Given any Hermitian manifold $(M^n,g)$ and under any local unitary frame $e$, the components of Bismut curvature and Chern curvature are related by
\begin{equation*} \label{eq:curvature2}
R^b_{i\bar{j}k\bar{\ell}} -  R_{i\bar{j}k\bar{\ell}} =  (T^{\ell}_{ik,\bar{j}} + \overline{   T^{k}_{j\ell ,\bar{i}}   }) + \sum_r \{ - T^r_{ik} \overline{T^r_{j\ell } }   +  T^{\ell}_{ir} \overline{T^k_{jr} }  - T^{j}_{ir} \overline{T^k_{\ell r} } -  T^{\ell}_{kr} \overline{T^i_{jr} } \},  
\end{equation*}
where the subscripts after comma stand for covariant derivatives with respect to the Bismut connection $\nabla^b$.
\end{lemma}

Taking the difference between the formula in Lemma \ref{lemma1} and Lemma \ref{lemma2}, we get the following formula which will be used later in our proof of Conjecture \ref{conj2} for BTP manifolds:

\begin{lemma} \label{lemma3}
Given any Hermitian manifold $(M^n,g)$ and under any local unitary frame $e$, the components of Riemannian curvature and Bismut curvature are related by
\begin{equation*} \label{eq:curvature2}
R^r_{i\bar{j}k\bar{\ell}} -  R^b_{i\bar{j}k\bar{\ell}} =  -\frac{1}{2} (T^{\ell}_{ik,\bar{j}} + \overline{   T^{k}_{j\ell ,\bar{i}}   }) + \frac{1}{4} \sum_r \{ - T^r_{ik} \overline{T^r_{j\ell } }   -  T^{\ell}_{ir} \overline{T^k_{jr} } -  T^{j}_{kr} \overline{T^i_{\ell r} }  +2 T^{j}_{ir} \overline{T^k_{\ell r} } +2 T^{\ell}_{kr} \overline{T^i_{jr} } \} , 
\end{equation*}
where subscripts after comma stand for covariant derivatives with respect to the Bismut connection $\nabla^b$.
\end{lemma}

Next let us recall the family of Gauduchon connections which is the line joining the Chern and Bismut connection:
$$ \nabla^{(t)}=(1-t)\nabla + t\nabla^b, \ \ \ \ \ t \in {\mathbb R}.  $$ 
It will be called the $t$-Gauduchon connection, and its curvature will be denoted by $R^{(t)}$. A beautiful recent result by Lafuente and Stanfield \cite{LS} states that any compact $t$-Gauduchon flat manifold must be K\"ahler, provided that $t\neq 0,1$. In fact, what they proved is a stronger statement, namely, for any compact Hermitian manifold and any $t \neq 0,1$,  the metric must be K\"ahler if it is assumed to be $t$-Gauduchon K\"ahler-like (in the sense that $R^{(t)}$ obeys all K\"ahler symmetries). 

For our later use, we will need the relation between components of $t$-Gauduchon curvature and Bismut curvature below (note that the $t=0$ case is just Lemma \ref{lemma2}). The proof is a straight-forward computation but we include it here for readers' convenience:

\begin{lemma} \label{lemma4}
Given any Hermitian manifold $(M^n,g)$ and under any local unitary frame $e$, the components of the $t$-Gauduchon curvature and Bismut curvature are related by
\begin{equation*} \label{eq:curvature2}
R^{(t)}_{i\bar{j}k\bar{\ell}} -  R^b_{i\bar{j}k\bar{\ell}} = (t\!-\!1) (T^{\ell}_{ik,\bar{j}} + \overline{   T^{k}_{j\ell ,\bar{i}}   })- (t\!-\!1) \sum_r \{ T^{j}_{ir} \overline{T^k_{\ell r} } + T^{\ell}_{kr} \overline{T^i_{jr} } \}
+ (t\!-\!1)^2 \sum_r \{  T^r_{ik} \overline{T^r_{j\ell } }   -  T^{\ell}_{ir} \overline{T^k_{jr} } \} , 
\end{equation*}
where subscripts after comma stand for covariant derivatives with respect to the Bismut connection $\nabla^b$.
\end{lemma}

\begin{proof} 
Let us fix a point $p\in M$. In a neighborhood of $p$, choose a local unitary frame $e$ so that the matrix of Bismut connection $\theta^b|_p=0$. This is always possible by \cite[Lemma 4]{YangZ}, where Chern connection was used, but the same proof works for any Hermitian connection obviously.  Denote by $\varphi$ the coframe dual to $e$. At the point $p$, the connection matrix  $\theta$ for Chern connection equals to $\theta|_p=\theta^b|_p - \gamma|_p= -\gamma|_p$, so by the structure equation $d\varphi =-\,^t\!\theta \wedge \varphi + \tau$, where $\varphi$ is understood to be a column vector and $\tau_r=\frac{1}{2}\sum_{i,k} T^r_{ik}\varphi_i \wedge \varphi_k$ is the column vector of Chern torsion, we get
\begin{equation*} \label{eq:partialbarvarphi}
\overline{\partial} \varphi_r = - \sum_{i,j} \overline{T^i_{jr}} \, \varphi_i \wedge \overline{\varphi}_j \ \ \mbox{at} \ p.
\end{equation*}
Denote by $\theta^{(t)}$, $\Theta^{(t)}$ the matrices of connection and curvature for the $t$-Gauduchon connection $\nabla^{(t)}$. We have $\theta^{(t)}= \theta^b + (t\!-\!1)\gamma$ and $\Theta^{(t)}= d\theta^{(t)} - \theta^{(t)}\wedge \theta^{(t)}$. So at the point $p$ we have
$$ \Theta^{(t)}-\Theta^b = (t\!-\!1) d\gamma -(t\!-\!1)^2\gamma \wedge \gamma .$$ 
Denote by $\gamma'$ the $(1,0)$-part of $\gamma$. We have $\gamma = \gamma' - \gamma'^{\ast}$ where $\ast$ means conjugate transpose. Taking the $(1,1)$-part in the above equation, we get
\begin{equation} \label{eq:tb}
(\Theta^{(t)}-\Theta^b)^{1,1} = (t\!-\!1) (\overline{\partial} \gamma' - \partial \gamma'^{\ast} ) +(t\!-\!1)^2(\gamma' \wedge \gamma'^{\ast}  + \gamma'^{\ast} \wedge \gamma') .
\end{equation} 
We compute at $p$ that
$$ \overline{\partial} \gamma'_{k\ell} = \overline{\partial}(T^{\ell}_{kr}\varphi_r) =  T^{\ell}_{kr,\bar{j}} \overline{\varphi}_j \wedge \varphi_r - T^{\ell}_{kr} \overline{  T^i_{jr}  } \varphi_i \wedge \overline{\varphi}_j  ,$$
where subscript after comma stands for covariant derivative with respect to the Bismut connection $\nabla^b$. Plugging this in (\ref{eq:tb}) we get the desired equation in Lemma \ref{lemma4}.
\end{proof}

We conclude this section by recalling a couple of properties for BTP manifolds. The first one is about the behaviour of Bismut curvature $R^b$ for BTP manifold (see \cite{ZhaoZ22}):

\begin{lemma} \label{lemma5} 
Let $(M^n,g)$ be a BTP manifold. Then under any local unitary frame $e$ one always has $R^b_{ijk\bar{\ell}}=0$, $R^b_{i\bar{j}k\bar{\ell}}=R^b_{k\bar{\ell}i\bar{j}}$, and
$ Q_{i\bar{j}k\bar{\ell}}: =R^b_{i\bar{j}k\bar{\ell}}- R^b_{k\bar{j}i\bar{\ell}} $  is always equal to
$$ Q=  \sum_r \{ - T^r_{ik} \overline{T^r_{j\ell } }   - T^{j}_{ir} \overline{T^k_{\ell r} } -  T^{\ell}_{kr} \overline{T^i_{jr} } +  T^{\ell}_{ir} \overline{T^k_{jr} } + T^{j}_{kr} \overline{T^i_{\ell r} }    \}  .$$
\end{lemma}

The second one is about the existence of {\em admissible frames} on non-balanced BTP manifolds \cite[Cor 5.2]{ZhaoZ22}:

\begin{lemma} \label{lemma6} 
Let $(M^n,g)$ be any non-balanced BTP manifold. Then locally there always exist admissible frames, namely, a local unitary frame $e$ with dual coframe $\varphi$ so that 
$$\eta = \lambda \varphi_n, \ \ \ T^n_{ik}=0, \ \ \  T^j_{in}=a_i\delta_{ij},\ \ \ \ R^b_{i \bar{j} k \bar{n}}=0 , \ \ \ \ \ \forall \ 1\leq i,j,k\leq n, $$
with constants $\lambda >0$, $a_n=0$,  and $a_1, \ldots , a_{n-1}$ satisfying  $a_1 + \cdots + a_{n-1}=\lambda$.
\end{lemma}

Here $\eta$ is Gauduchon's torsion $1$-form (\cite{Gauduchon}), which is the global $(1,0)$-form on the manifold defined by $d(\omega^{n-1})=-\eta \wedge \omega^{n-1}$ where $\omega$ is the K\"ahler form of $g$. By definition, $g$ is balanced if and only if $\eta =0$. So when the metric is non-balanced BTP, we have $\nabla^be_n=0$, thus the condition $R^b_{\ast \bar{\ast} \ast \bar{n}}=0$ holds.

\vspace{0.3cm}

\section{BTP manifolds with constant holomorphic sectional curvature}

In this section, let us assume that $(M^n,g)$ is a BTP manifold, with constant holomorphic sectional curvature either in Riemannian connection $\nabla^r$ or in $t$-Gauduchon connection $\nabla^{(t)}$. Note that the $t=0$ case is just the Chern connection.

For a tensor $P_{i\bar{j}k\bar{\ell}}$, let us define its symmetrization by
$$ \widehat{P}_{i\bar{j}k\bar{\ell}} = \frac{1}{4}\{ P_{i\bar{j}k\bar{\ell}} + P_{k\bar{j}i\bar{\ell}} + P_{i\bar{\ell}k\bar{j}} + P_{k\bar{\ell}i\bar{j}}\} .$$
Clearly, the condition $P_{X\overline{X}X\overline{X}}=c|X|^4$ for any $X$ is equivalent to
$$  \widehat{P}_{i\bar{j}k\bar{\ell}}  = \frac{c}{2} ( \delta_{ij}\delta_{k\ell} + \delta_{i\ell} \delta_{kj} ) , \ \ \ \ \ \forall \ 1\leq i,j,k,\ell \leq n. $$
Therefore we have
\begin{eqnarray*}
H^r=c  \ \ \, & \Longleftrightarrow & \ \ \widehat{R^r}_{i\bar{j}k\bar{\ell}} = \frac{c}{2} ( \delta_{ij}\delta_{k\ell} + \delta_{i\ell} \delta_{kj} ), \\
H^{(t)}=c   \ & \Longleftrightarrow & \ \widehat{R^{(t)}}_{i\bar{j}k\bar{\ell}} = \frac{c}{2} ( \delta_{ij}\delta_{k\ell} + \delta_{i\ell} \delta_{kj} ), 
\end{eqnarray*}
for any $1\leq i,j,k,\ell \leq n$. To simplify the writing, let us introduce the following short-hand notations:
\begin{equation*}
w=  T^r_{ik} \overline{T^r_{j\ell } } , \ \ \ v^j_i =  T^j_{ir} \overline{T^k_{\ell r} } , \ \ \  v^{\ell}_k = T^{\ell}_{kr} \overline{T^i_{jr } }, \ \ \   v^{\ell}_i =  T^{\ell}_{ir} \overline{T^k_{jr } } , \ \ \ v^j_k =   T^j_{kr} \overline{T^i_{\ell r} } ,
\end{equation*}
where $1\leq i,j,k,\ell \leq n$ and $r$ is summed from $1$ to $n$. If we perform the symmetrization to these terms, we get
\begin{equation*}
\widehat{^{\,}w}=0,  \ \ \ \ \ \ \widehat{v^j_i}= \widehat{v^{\ell}_k} =\widehat{v^{\ell}_i} = \widehat{v^j_k} = \frac{1}{4} (v^j_i +  v^{\ell}_k + v^{\ell}_i + v^j_k  ).
\end{equation*}

From this we deduce the following:

\begin{lemma} \label{lemma7}
Suppose that $(M^n,g)$ is a BTP manifold with constant Riemannian holomorphic sectional curvature $H^r=c$. Then under any local unitary frame  $e$ it holds
$$ R^b_{   i\bar{j} k\bar{\ell}   } = \frac{c}{2} ( \delta_{ij}\delta_{k\ell} + \delta_{i\ell} \delta_{kj} )  -\frac{1}{2}w -\frac{5}{8}(v^j_i+v^{\ell}_k) + \frac{3}{8}(v^{\ell}_i + v^j_k). $$
\end{lemma}

\begin{proof}
For BTP manifolds, by Lemma \ref{lemma5} we know that the Bismut curvature tensor $R^b$ always satisfies the symmetry condition $R^b_{i\bar{j}k\bar{\ell}} = R^b_{k\bar{\ell}i\bar{j}}$ and 
$$ Q_{i\bar{j}k\bar{\ell}} : = \, R^b_{i\bar{j}k\bar{\ell}}- R^b_{k\bar{j}i\bar{\ell}} \, = \, - w - v^j_i -  v^{\ell}_k + v^{\ell}_i + v^j_k .$$
Therefore for BTP manifolds it holds that
\begin{equation} \label{eq:Rbhat}
 \widehat{ R^b}_{i\bar{j}k\bar{\ell}} = \frac{1}{2}(R^b_{i\bar{j}k\bar{\ell}}+ R^b_{k\bar{j}i\bar{\ell}}) = \frac{1}{2}( 2R^b_{i\bar{j}k\bar{\ell}} - Q_{i\bar{j}k\bar{\ell}}) = R^b_{i\bar{j}k\bar{\ell}} + \frac{1}{2}(w + v^j_i +  v^{\ell}_k - v^{\ell}_i - v^j_k). 
\end{equation}
Under the assumption $H^r=c$, we have 
$  \widehat{R^r}_{i\bar{j}k\bar{\ell} } = \frac{c}{2} ( \delta_{ij}\delta_{k\ell} + \delta_{i\ell } \delta_{kj} )$. On the other hand, by taking the symmetrization on the formula in Lemma \ref{lemma3} we get 
$$ \widehat{R^r}- \widehat{R^b} = \frac{1}{8}(v^j_i +  v^{\ell}_k + v^{\ell}_i + v^j_k).$$
Now plugging in the result for   $\widehat{R^r}$ and using (\ref{eq:Rbhat}), we obtain the desired expression for $R^b$ stated in the lemma.
\end{proof}

Similarly, by taking the symmetrization of the formula in Lemma \ref{lemma4} we get
$$ \widehat{R^{(t)}} -\widehat{ R^b} = \frac{1}{4}(1-t^2)(v^j_i +  v^{\ell}_k + v^{\ell}_i + v^j_k). $$
Combine this will (\ref{eq:Rbhat}) we obtain the expression of $R^b$ in $\widehat{R^{(t)}} $ and torsion terms, and get
\begin{lemma} \label{lemma8}
Suppose that $(M^n,g)$ is a BTP manifold with constant $t$-Gauduchon holomorphic sectional curvature $H^{(t)}=c$. Then under any local unitary frame  $e$ it holds
$$ R^b_{   i\bar{j} k\bar{\ell}   } = \frac{c}{2} ( \delta_{ij}\delta_{k\ell} + \delta_{i\ell} \delta_{kj} )  -\frac{1}{2}w +\frac{t^2-3}{4}(v^j_i+v^{\ell}_k) + \frac{t^2+1}{4}(v^{\ell}_i + v^j_k). $$
\end{lemma}

Note that the  curvature tensor $R_{i\bar{j}k\bar{\ell}}$ of any K\"ahler metric is symmetric with respect to its first and third position, as well as  its second and fourth position. As a result we have $\widehat{R}=R$ and the holomorphic sectional curvature $H$ would determine the entire $R$ algebraically, and we ended up with the three complex space forms when $H$ is constant. In the Hermitian case, however, the curvature tensor of Riemannian or $t$-Gauduchon connection no longer obey those K\"ahler symmetries in general, and holomorphic sectional curvature $H$ could only determine the symmetrization part of $R$ instead. This is the main difficulty underlying the Holomorphic Sectional Curvature Conjecture.

Nonetheless, the above two lemmas tell us that for BTP manifolds, when either the Riemannian connection or a $t$-Gauduchon connection  has constant holomorphic sectional curvature, then its entire Bismut curvature $R^b$ (hence other curvatures as well) can be expressed in terms of the Chern torsion components along. With these formula in hand, we are now ready to prove Theorems \ref{thm1} and \ref{thm4} in the non-balanced case.

\begin{proof}[{\bf Proof of Theorems \ref{thm1} and \ref{thm4} (non-balanced case).}]
Assume that $(M^n,g)$ is a compact, non-balanced BTP manifold. By Lemma \ref{lemma6}, there always exist admissible frames. Let us assume that the unitary frame $e$ is an admissible one. Hence 
$$T^n_{ik}=0, \ \ \ T^j_{in}=\delta_{ij}a_i, \ \ \ R^b_{i\bar{j}k\bar{n}}=0, \ \ \ \ \forall \ 1\leq i,j,k\leq n,  $$
where $a_i$ are constants satisfying $a_n=0$ and  $a_1+\cdots +a_{n-1}=\lambda >0$. 

Assume now that the Riemannian holomorphic sectional curvature $H^r=c$ is a constant. By letting $i=j\leq k=\ell =n$ in the formula in Lemma \ref{lemma7}, we get 
$$ 0 = \frac{c}{2}(1+\delta_{in}) + (- \frac{1}{2} +\frac{3}{8}) |a_i|^2. $$
Letting $i=n$ we deduce $c=0$, and then $a_i=0$ for each $i$, contradicting with the fact that their sum is equal to $\lambda >0$. This shows that for any non-balanced BTP manifolds, its $H^r$ cannot be a constant. 

Similarly, if we assume that the $t$-Gauduchon connection has constant holomorphic sectional curvature $H^{(t)}=c$. Then by choosing $i=j\leq k=\ell =n$ in the formula in Lemma \ref{lemma8}, we get 
$$ 0 = \frac{c}{2}(1+\delta_{in}) + (- \frac{1}{2} +\frac{t^2+1}{4}) |a_i|^2. $$
Letting $i=n$ we deduce $c=0$, therefore $0=\frac{t^2-1}{4} |a_i|^2$. If $t\neq \pm 1$, then again it leads to $a_i=0$ for all $i$ and we get a contradiction. So for any non-balanced BTP manifolds and for any $t\neq \pm 1$, the $t$-Gauduchon holomorphic sectional curvature $H^{(t)}$ cannot be a constant. Note that when  $t=0$ it is just the Chern connection case. For the $t=1$ or $-1$ cases, by letting $i=j=k=\ell =n$ in the formula in either Lemma \ref{lemma7} or Lemma  \ref{lemma8}, we end up with $c=0$. So for non-balanced BTP manifolds, the holomorphic sectional curvature of $\nabla^b$ or $\nabla^{(-1)}$ cannot be a non-zero constant. This completes the proof of  Theorems \ref{thm1} and \ref{thm4} in the case when the BTP manifold is assumed to be non-balanced.
\end{proof}
As we mentioned before, when $n\geq 3$, BTP manifolds can be balanced (and non-K\"ahler). Balanced BTP manifolds form a highly restrictive and interesting class. In this case, the Bismut curvature tensor does not always have a non-trivial kernel, and the formula from Lemma \ref{lemma7} or \ref{lemma8} no longer leads us to a quick conclusion, so we do not know how to prove the conjectures at the present time, except in the special case when the dimension is $3$. In this case, the recent work of Zhao and the second named author \cite{ZhaoZ22} provided a classification for all balanced BTP threefolds, so we could just verify the conjectures in a  case by case manner. In the following we will first recall this classification result for balanced BTP threefolds.

Let $(M^3,g)$ be a  balanced, non-K\"ahler compact BTP Hermitian threefold. By the observation in \cite{ZhouZ} and \cite{ZhaoZ22}, for any given point $p\in M$, there always exists a unitary frame $e$ (which will be called {\em special frames} from now on) in a neighborhood of $p$, so that under $e$ the only possibly non-zero Chern torsion components are $a_i=T^i_{jk}$, where $(ijk)$ is a cyclic permutation of $(123)$. Furthermore, each $a_i$ is a global constant on $M^3$, with $a_1=\cdots =a_r>0$, $a_{r+1}=\cdots =0$, where $r\in \{ 1,2,3\}$   is the rank of the $B$ tensor, which is the global $2$-tensor on any Hermitian manifold defined under any unitary frame by
$ B_{i\bar{j}} = \sum_{k,\ell } T^j_{k\ell} \overline{   T^i_{k\ell }  }$.
The conclusion in \cite{ZhaoZ22} indicates that any compact balanced (but non-K\"ahler) BTP threefold must be one of the following:

\begin{itemize}
\item $r=3$, $(M^3,g)$ is a compact quotient of the complex simple Lie group $SO(3,{\mathbb C})$, in particular it is Chern flat;

\item  $r=1$, $(M^3,g)$ is the so-called Wallach threefold, namely, $M^3$ is biholomorphic to the flag variety ${\mathbb P}(T_{{\mathbb P}^2} )$ while $g$ is the K\"ahler-Einstein metric minus the square of the null-correlation section. Scale $g$ by a positive constant if necessary, the Bismut curvature matrix under a special frame $e$ is
\begin{equation} \label{eq:Fanotype}
\Theta^b = \left[ \begin{array}{ccc} \alpha+\beta & 0 & 0\\ 0 & \alpha & \sigma \\ 0 & - \bar{\sigma} & \beta \end{array} \right],  \ \ \ \ \ \ \ \ \ \left\{ \begin{array}{lll}  \alpha = \varphi_{1\bar{1}}+ (1\!-\!b)\varphi_{2\bar{2}} + b \varphi_{3\bar{3}} +t \varphi_{2\bar{3}} + \bar{t} \varphi_{3\bar{2}},
 && \\ \beta = \varphi_{1\bar{1}}+ b\varphi_{2\bar{2}} + (1\!-\!b) \varphi_{3\bar{3}} -t \varphi_{2\bar{3}} - \bar{t} \varphi_{3\bar{2}},  & & \\ \sigma = t\varphi_{2\bar{2}} -t \varphi_{3\bar{3}} +s \varphi_{2\bar{3}} + (1\!+\!b) \varphi_{3\bar{2}}, & &
\end{array} \right.
\end{equation}
where $b$ is a real constant, $t,s$ are complex constants, $\varphi$ is the coframe dual to $e$, and we wrote $\varphi_{i\bar{j}}$ for $\varphi_i \wedge \overline{\varphi}_j$. 

\item $r=2$, in this case $(M^3,g)$ is said to be of {\em middle type}. Again under appropriate scaling of the metric, the Bismut curvature matrix under $e$ becomes
\begin{equation} \label{eq:middletype}
\Theta^b = \left[ \begin{array}{ccc} d\alpha & d\beta_0 & \\ - d\beta_0 & d\alpha & \\ & & 0 \end{array} \right], \ \ \ \ \ \ \ \ \ \left\{ \begin{array}{ll}  d\alpha = \ x(\varphi_{1\bar{1}}+\varphi_{2\bar{2}}) +iy(\varphi_{2\bar{1}}- \varphi_{1\bar{2}}), & \\ d\beta_0 = -iy (\varphi_{1\bar{1}}+\varphi_{2\bar{2}}) + (x\!-\!2)(\varphi_{2\bar{1}}-\varphi_{1\bar{2}}), & \end{array} \right. 
\end{equation}
where $x,y$ are real-valued local smooth functions. 
\end{itemize}

With these explicit information on Bismut curvature at hand, we are now ready to finish the proof of Theorems \ref{thm1} and \ref{thm4}.

\begin{proof}[{\bf Proof of Theorems \ref{thm1} and \ref{thm4} (balanced threefold case).}]
Let $(M^3,g)$ be a compact balanced, non-K\"ahler BTP threefold. Then it must be one of three types above according to its $B$ rank $r$. Here let us scale the metric $g$ so that in the $r=1$ and $r=2$ case the Bismut curvature components under a special frame $e$ are given by formula (\ref{eq:Fanotype}) and (\ref{eq:middletype}), respectively.

First let us assume that its Riemannian connection has constant holomorphic sectional curvature: $H^r=c$. Then $R^b$ is given by Lemma \ref{lemma7}. By letting $i=j$ and $k=\ell$, we obtain
\begin{equation*} 
R^b_{i\bar{i}k\bar{k}} = \frac{c}{2}(1+\delta_{ik}) +\sum_s \{ -\frac{1}{2}  |T^s_{ik}|^2  + \frac{3}{8}( |T^k_{is}|^2 + |T^i_{ks}|^2 )\}.
\end{equation*}
Note that here we used the fact that $T^i_{i\ast }=0$ since $e$ is special. When $i\neq k$, the sigma part in the right hand side of the above equation becomes $-\frac{1}{2}a_j^2 + \frac{3}{8}(a_k^2+ a_i^2)$, where $j\neq i,k$. Therefore we have
\begin{equation} \label{eq:Hr1}
\left\{ \begin{array}{ll} R^b_{i\bar{i}i\bar{i}}=c, & \\
 R^b_{i\bar{i}k\bar{k}} = \frac{c}{2}-\frac{1}{2}a_j^2 + \frac{3}{8}(a_k^2+ a_i^2),  \ \ \ \mbox{if} \ \{ i,j,k\} =\{ 1,2,3\} & \end{array}\right. 
\end{equation}
Now if the $B$ rank of $M^3$ is $r=3$, then $a_1=a_2=a_3>0$ and $g$ is Chern flat. By Lemma \ref{lemma2} we have
$R^b = -w + v^{\ell}_i - v^j_i - v^{\ell}_k$, therefore
$$ R^b_{i\bar{i}k\bar{k}} = \sum_s \{ -|T^s_{ik}|^2  + |T^k_{is}|^2 \}. $$
When $i=k$, it is clearly equal to zero, while when $i\neq k$, it equals to $-|a_j|^2 + |a_k|^2=0$, where $\{i,j,k\} =\{ 1,2,3\}$. That is, $R^b_{i\bar{i}k\bar{k}}=0$ for all $1\leq i,k\leq 3$. Compare with (\ref{eq:Hr1}), we get $c=0$  and $0=\frac{1}{4}a_1^2$, which is a contradiction.

If the $B$ rank of $M^3$ is $1$, then $a_1>0=a_2=a_3$, and by (\ref{eq:Hr1}) we have $R^b_{1\bar{1}1\bar{1}}=c$ and $R^b_{2\bar{2}1\bar{1}}=\frac{c}{2} +\frac{3}{8}a_1^2$. On the ther hand, by (\ref{eq:Fanotype}) we know that
$$ \Theta^b_{11} = \alpha +\beta = 2\varphi_{1\bar{1}}+ \varphi_{2\bar{2}}+\varphi_{3\bar{3}},$$ 
therefore $R^b_{1\bar{1}1\bar{1}}=2R^b_{2\bar{2}1\bar{1}}$, or
$$ c =2(\frac{c}{2} +\frac{3}{8}a_1^2),$$
which leads to $a_1=0$, a contradiction.

If the $B$ rank of $M^3$ is $2$, then we are in the middle type with $a_1=a_2>0=a_3$. By (\ref{eq:Hr1}) we have $R^b_{3\bar{3}3\bar{3}}=c$ and 
$$R^b_{2\bar{2}3\bar{3}}=\frac{c}{2} -\frac{1}{2}a_1^2 +\frac{3}{8}a_2^2 = \frac{c}{2} -\frac{1}{8}a_1^2.$$ 
On the other hand, by (\ref{eq:middletype}), $\Theta^b_{33}=0$, so $R^b_{2\bar{2}3\bar{3}}=R^b_{3\bar{3}3\bar{3}}=0$. Hence $c=0$ and $a_1=0$, again a contradiction.

The above discussion shows that, under the assumption that $H^r=c$, we will have the expression (\ref{eq:Hr1}) for $R^b$, which will contradict with each of the three cases in the classification of balanced BTP threefolds. Thus for any compact balanced (non-K\"ahler) BTP threefolds, its Riemannian connection cannot have constant holomorphic sectional curvature. 

Next, let us assume that $t\in {\mathbb R}$ and the $t$-Gauduchon connection has constant holomorphic sectional curvature: $H^{(t)}=c$. Under a special frame $e$, the only possibly non-zero components of Chern torsion are $T^1_{23}=a_1$, $T^2_{31}=a_2$, and $T^3_{12}=a_3$. By letting $i=j$ and $k=\ell$ in the formula in Lemma \ref{lemma8}, we get 
\begin{equation} \label{eq:Ht}
\left\{ \begin{array}{ll} R^b_{i\bar{i}i\bar{i}}=c, & \\
 R^b_{i\bar{i}k\bar{k}} = \frac{c}{2}-\frac{1}{2}a_j^2 + \frac{t^2+1}{4}(a_k^2+ a_i^2),  \ \ \ \mbox{if} \ \{ i,j,k\} =\{ 1,2,3\} .& \end{array}\right. 
\end{equation}

If the $B$ rank of $M^3$ is $r=3$, then $a_1=a_2=a_3>0$ and $g$ is Chern flat. As before, we deduce from Lemma \ref{lemma2} that $R^b_{i\bar{i}k\bar{k}} =0$. On the other hand, the equation (\ref{eq:Ht}) becomes $R^b_{i\bar{i}i\bar{i}} = c$ and $R^b_{i\bar{i}k\bar{k}} = \frac{c}{2} + \frac{t^2}{2}a_1^2$ when $i\neq k$. Thus we get $c=0$ and $t=0$. Therefore, for any $t\neq 0$, the $t$-Gauduchon connection cannot have constant holomorphic sectional curvature in the $r=3$ case. (While for $t=0$, namely for the Chern connection, the manifold is already Chern flat so Conjecture \ref{conj1} holds). 

If the $B$ rank of $M^3$ is $r=1$, then $a_1>0=a_2=a_3$ and $g$ is the Wallach threefold. As we have seen before, in this case $R^b_{1\bar{1}1\bar{1}}= 2R^b_{2\bar{2}1\bar{1}}$. But (\ref{eq:Ht}) leads to $R^b_{1\bar{1}1\bar{1}}=c$ and $R^b_{2\bar{2}1\bar{1}}= \frac{c}{2} + \frac{t^2+1}{4}a_1^2$. Therefore we get  $a_1=0$, a contradiction. So for the Wallach threefold, $H^{(t)}$ cannot be a constant for any $t$. 

If the $B$ rank of $M^3$ is $r=2$, then we are in the middle type and $a_1=a_2>0=a_3$. In this case Equation (\ref{eq:Ht}) gives us $R^b_{3\bar{3}3\bar{3}}=c$ and 
$$ R^b_{2\bar{2}3\bar{3}}= \frac{c}{2} - \frac{1}{2}a_1^2 + \frac{t^2+1}{4}a_2^2 =  \frac{c}{2} +\frac{t^2-1}{4}a_1^2.$$
On the hand, by (\ref{eq:middletype}) we have $\Theta^b_{33}=0$, so $c=0$ and 
$ (t^2-1)a_1^2=0$, which leads to a contradiction when $t\neq \pm 1$. 

For $t=1$, namely the Bismut connection case, by (\ref{eq:middletype}) we have $R^b_{3\bar{3}3\bar{3}}=0$  and  $R^b_{1\bar{1}1\bar{1}}=x$. So if $H^b$ is a constant, then the constant must be $0$ and $x=0$. Given any tangent direction $X=\sum_i x_ie_i$, write $p=|x_1|^2+|x_2|^2$, $q=x_1\overline{x}_2 - x_2\overline{x}_1$. By (\ref{eq:middletype}) we have
$$ R^b_{X\bar{X}X\bar{X}} = p \,d\alpha (X, \bar{X}) + q \,d\beta_0(X,\bar{X}) = p(-iyq)+q(-iyp+2q) = 2q \big( q-iyp). $$ 
Clearly, for any fixed $y$ value, this expression cannot be identically zero for all choices of $X$. That is,  $H^b$ can never be a constant. Similarly, for the `mirror' connection $\nabla^{(-1)}$, by letting $t=-1$ in Lemma \ref{lemma4} we get
$$ R^{(-1)} - R^b = 2(v^j_i + v^{\ell}_k) +4 (w - v^{\ell}_i), $$
for any BTP manifold. So for any balanced BTP threefold and under any special frame, we have $ R^{(-1)}_{i\bar{i}i\bar{i}}=R^b_{i\bar{i}i\bar{i}}$, thus in the middle type case $H^{(-1)}$ again cannot be constant. 

In summary, we have proved that, for any compact balanced, non-K\"ahler BTP threefold $(M^3,g)$, its Riemannian connection or $t$-Gauduchon connection cannot have constant holomorphic sectional curvature, with the exception of $t=0$ (Chern connection) and the $B$ rank $r=3$ case. But in this case the manifold is already Chern flat. Thus we have completed the proof of Theorems \ref{thm1} and \ref{thm4}. 
\end{proof}

As we have noted in the proofs, for BTP manifold, the $t=-1$ case of Lemma \ref{lemma4} gives us
$$ R^{(-1)} - R^b = 2(v^j_i + v^{\ell}_k) + 4(w- v^{\ell}_i). $$
Hence $\widehat{ R^{(-1)}} - \widehat{R^b}=0$. So for BTP manifolds, the condition $H^b=0$ is equivalent to $H^{(-1)}=0$. In particular, isosceles Hopf manifolds have $H^{(-1)}=0$. On the other hand, $ R^{(-1)} \not\equiv 0$ for such manifolds by Lafuente-Standfield Theorem (or by direct verification). So in Conjecture \ref{conj4} one has to exclude the $\nabla^{(-1)}$ in addition to $\nabla^{(1)}=\nabla^b$. 

It would be an interesting question to classify all (non-balanced) BTP manifolds  with $H^b=0$ (which is equivalent to $H^{(-1)}=0$), other than the Bismut flat ones or the isosceles Hopf manifolds.

\vs

\noindent\textbf{Acknowledgments.} The second named author would like to thank Haojie Chen, Xiaolan Nie, Kai Tang, Bo Yang, Yashan Zhang, and Quanting Zhao for their interest and helpful discussion.

\vs

\end{document}